\documentclass[12pt]{article}
\usepackage{amssymb,amsfonts,amsmath, psfrag,eepic,colordvi}
\topskip  
\numberwithin{equation}{section}
\newtheorem{theorem}{Theorem}[section]

\newtheorem{conjecture}[theorem]{Conjecture}

\newtheorem{lemma}[theorem]{Lemma}

\newtheorem{example}[theorem]{Example}

\begin{document}
\parskip 6pt

\pagenumbering{arabic}
\def\sof{\hfill\rule{2mm}{2mm}}
\def\ls{\leq}
\def\gs{\geq}
\def\SS{\mathcal S}
\def\qq{{\bold q}}
\def\MM{\mathcal M}
\def\TT{\mathcal T}
\def\EE{\mathcal E}
\def\lsp{\mbox{lsp}}
\def\rsp{\mbox{rsp}}
\def\pf{\noindent {\it Proof.} }
\def\mp{\mbox{pyramid}}
\def\mb{\mbox{block}}
\def\mc{\mbox{cross}}
\def\qed{\hfill \rule{4pt}{7pt}}
\def\block{\hfill \rule{5pt}{5pt}}

\begin{center}
{\Large\bf  On a conjecture about enumerating $(2+2)$-free posets  }
\vskip 6mm
\end{center}

\begin{center}
{\small   Sherry H. F. Yan \\[2mm]
 Department of Mathematics, Zhejiang Normal University, Jinhua
321004, P.R. China
\\[2mm]
 huifangyan@hotmail.com
  \\[0pt]
}
\end{center}

\noindent {\bf Abstract.} Recently, Kitaev and Remmel posed
 a conjecture concerning  the   generating function for the number
 of unlabeled
$(2+2)$-free posets  with respect to number of elements  and number
of minimal elements.  In this paper, we   present    a combinatorial
proof of  this conjecture.




\section{Introduction}
A  poset is said to be $(2+2)$-free if it does not contain an
induced subposet that is isomorphic to $2+2$, the union of two
disjoint $2$-element chains. In a poset, let $D(x)$ be the set of
{\em predecessors} of  an element $x$ (the strict down-set of $x$).
Formally, $D(x)=\{y: y<x\}$.  A poset $P$  is $(2+2)$-free if and
only if its sets of predecessors, $D(P)=\{D(x): x\in P\}$ can be
 written as
  $$
  D(P)=\{D_0, D_1, \ldots, D_k \}
 $$
where $\emptyset=D_0\subset D_1\subset\ldots \subset D_k$, see
\cite{Bog, melon}. In such context, we say that $x\in P$ has {\em
level} $i$ if $D(x)=D_i$.  An element $x$ is said to be a {\em
minimal} element if $x$ has level $0$.

Let $p_n$ be the number of unlabeled  $(2+2)$-free posets on $n$
elements. EI-Zahar \cite{zahar} and Khamis \cite{khamis} used a
recursive description of  $(2+2)$-free posets to derive a pair of
functional equations that define the generating function for the
number $p_n$. But they did not solve these equations. Recently,
using functional equations and the Kernel method,
Bousquet-M$\acute{e}$lou et al. \cite{melon} showed that the
generating function for the number $p_n$ of unlabeled  $(2+2)$-free
posets on n elements is given by
\begin{equation} \label{pt}
P(t)= \sum_{n\geq 0}p_{n} t^n =  \sum_{n\geq
0}\prod_{i=1}^{n}(1-(1-t)^{i} ).
\end{equation}
Note that throughout this paper, the empty product as usual is taken
to be $1$.  In fact, they studied a more general function of
unlabeled $(2+2)$-free posets according to number of elements,
number of levels and level of minimum maximal elements. Zagier
\cite{zag} proved that Formula (\ref{pt}) is also the generating
function for certain involutions introduced by Stoimenow
\cite{Stoi}.

Given a sequence of integers   $x=(x_1,x_2,\ldots, x_n)$, we say
that the sequence $x$ has an ascent at position $i$ if
$x_i<x_{i+1}$. The number of ascents of $x$ is denoted by $asc(x)$.
A sequence $x=(x_1,x_2, \ldots, x_n)$ is said to be an {\em ascent
sequence of length $n$} if it satisfies $x_1=0$ and $0\leq x_i\leq
asc(x_1,x_2,\ldots, x_{i-1})+1$ for all $2\leq i\leq n$. Ascent
sequences were introduce by Bousquet-M$\acute{e}$lou et al.
\cite{melon}  to unify three combinatorial structures.
Bousquet-M$\acute{e}$lou et al. \cite{melon} constructed bijections
between unlabeled  $(2+2)$-free posets and ascent sequences, between
ascent sequences and permutations avoiding a certain pattern,
between unlabeled $(2+2)$-free posets and a class of involutions
introduced by Stoimenow \cite{Stoi}.

Recently, Kitaev and Remmel \cite{Kitaev} extended the work of
Bousquet-M$\acute{e}$lou et al. \cite{melon}. They found generating
function for unlabeled  $(2+2)$-free posets when four statistics are
taken into account, one of which is the number of minimal elements
in a poset. The key strategy used by Bousquet-M$\acute{e}$lou et al.
\cite{melon} and Kitaev and Remmel \cite{Kitaev}  is to translate
statistics on $(2+2)$-free posets to statistics on ascent sequences
using  the bijection between unlabeled  $(2+2)$-free posets and
ascent sequences given by Bousquet-M$\acute{e}$lou et al.
\cite{melon}. Let $p_{n,k}$ be the number of $(2+2)$-free posets on
$n$ elements with $k$ minimal elements, with the assumption
$p_{0,0}=1$. Under the bijection between unlabeled $(2+2)$-free
posets and ascent sequences, the number of unlabeled $(2+2)$-free
posets on $n$ elements with $k$ minimal elements is equal to that of
ascent sequences of length $n$ with $k$ zeros. Kitaev and Remmel
\cite{Kitaev} derived that the generating function for the number
$p_{n,k}$ is given by
$$
P(t,z)=\sum_{n\geq 0, k\geq 0}p_{n,k}z^kt^n= 1+ \sum_{n\geq
0}{zt\over (1-tz)^{n+1}}\prod_{i=1}^{n}(1-(1-t)^i),
$$
by counting ascent sequences with respect to  length and number of
zeros. Moreover, they conjectured the function $P(t,z)$ can be
written in a simpler form.

\begin{conjecture}\label{conjecture}

\begin{equation}\label{identity}
P(t,z)=\sum_{n\geq 0, k\geq 0}p_{n,k}z^kt^n= \sum_{n\geq
0}\prod_{i=1}^{n}(1-(1-t)^{i-1}(1-zt)).\end{equation}
\end{conjecture}

   The objective of this paper  is to give a combinatorial proof
of Conjecture \ref{conjecture}.  In order to prove the conjecture,
we need  two more  combinatorial structures: upper triangular
matrices with non-negative  integer entries such that all rows and
columns contain at least one non-zero entry, which was introduced by
Dukes and Parviainen \cite{Duck},  and upper triangular
$(0,1)$-matrices in which all columns contain at least one non-zero
entry.

  Let
$\mathcal{A}_n$ be the collection of upper triangular matrices with
non-negative integer entries which sum to $n$.   A {\em
$(0,1)$-matrix} is a matrix in which each entry is either  $0$ or
$1$. Let $\mathcal{M}_n$ be the set of $(0,1)$-matrices in
$\mathcal{A}_n$ in which  all columns contain at least one non-zero
entry. Denote by $\mathcal{I}_n$ the set of matrices in
$\mathcal{A}_n$ in which all rows and columns contain at least one
non-zero entry.  Given a matrix $A$, denoted  by $A_{i,j}$ the entry
in row $i$ and column $j$.  Let $dim(A)$ be the number of rows in
the matrix $A$. The sum of all entries in row $i$ is called the {\em
row sum}  of row $i$, denoted by $rsum_i(A)$. The {\em column sum}
of column $i$, denoted by $csum_i(A)$,  can be defined similarly. A
row is said to be {\em zero} if its row sum is zero.

Let $A$ be a matrix in $\mathcal{M}_n$, define  $min_i(A)$  to be
the least   value of $j$ such $A_{j,i}$ is non-zero. A column $i$ of
$A$ is said to be {\em improper} if it satisfies one of the
following two cases: (1) $csum_i(A)\geq 2$; (2) for $1<i\leq
dim(A)$,  we have $csum_i(A)=1$, $rsum_i(A)=0$, and
$min_i(A)<min_{i-1}(A)$. Otherwise, column $i$ is said to be {\em
proper}. The matrix $A$ is said to be {\em improper} if there  is at
least one improper column in $A$; otherwise, the matrix $A$ is said
to be {\em proper}. Given an improper matrix $A\in \mathcal{M}_n$,
 define $index(A)$ to be the  largest  value   $i$ such that
 column
 $i$ is improper. Denote by $\mathcal{PM}_n$ the set of proper
 matrices in $\mathcal{M}_n$.

\begin{example}
Consider the following matrix $A\in \mathcal{M}_8$:
$$
A=\begin{bmatrix}
1 &0  &1  &0 &0& 1   \\
0 & 1 &0  & 1 &1& 0  \\
0 & 0& 0 & 0 &0 &0  \\
0 & 0 & 0& 0 & 0 &1  \\
0 & 0 & 0& 0 & 0& 0  \\
0 & 0& 0 &0 &0& 1  \\
\end{bmatrix}.
$$
We have $dim(A)=6, min_1(A)=1, min_2(A)=2, min_3(A)=1, min_4(A)=2,
min_5(A)=2, min_6(A)=1 $.  There are two improper columns, that is,
columns $3$  and  $6$.  Hence,  we have $index(A)=6$.
\end{example}

Denote by $\mathcal{PM}_{n,k}$ the set of matrices $A\in
\mathcal{PM}_n$ with $rsum_1(A)=k$ and $\mathcal{I}_{n,k}$ the set
of matrices $A\in \mathcal{I}_n$ with $rsum_1(A)=k$. Dukes and
Parviainen \cite{Duck} constructed a recursive  bijection between
the set $\mathcal{I}_n$ and the set of ascent sequences of length
$n$.    Under their bijection, they showed that the number of upper
triangular matrices $A\in \mathcal{I}_{n}$ with $rsum_1(A)=k$ is
equal to the number of ascent sequences of length $n$ with $k$
zeros, which implies that the cardinality of $\mathcal{I}_{n,k}$ is
also given by $p_{n,k}$. In this paper, we will prove Conjecture
\ref{conjecture} by showing that the generating function for the
number of matrices in $\mathcal{I}_{n,k}$ is given by the right-hand
side of Formula (\ref{identity}).

  In Section 2, we present a parity reversing and weight
preserving involution on the set  $\mathcal{M}_n\setminus
\mathcal{PM}_n$.   In Section 3, we  prove that the right-hand side
of Formula (\ref{identity}) is the generating function for the
number of matrices in $\mathcal{PM}_{n,k}$.   Moreover,  we show
that there is a bijection between the set $\mathcal{PM}_{n,k}$ and
the set $\mathcal{I}_{n,k}$ in answer to Conjecture
\ref{conjecture}.

\section{A parity reversing and weight preserving  involution}

In this section, we will construct a parity reversing and weight
preserving  involution on   the set $\mathcal{M}_n\setminus
\mathcal{PM}_n$. Before constructing the involution, we need some
definitions.

Given a matrix $A\in \mathcal{M}_n$, the {\em weight } of the matrix
 $A$ is assigned by
   $ z^{rsum_1(A)}$.
    Given  a subset $S$ of  the set $\mathcal{M }_n$, the {\em weight} of  $S$, denoted by $W(S)$,  is the  sum of the weights
 of all matrices in $S$.   We define the {\em parity} of the matrix $A$ to be
 the parity of the number $n-dim(A)$. Denote by $\mathcal{EM}_n$ (resp. $\mathcal{OM}_n$ ) the
 set of matrices in
 $\mathcal{M}_n$ whose parity are even (resp. odd).

 \begin{theorem}\label{involution}
There is a parity reversing  and weight preserving involution $\Phi$
on the set $\mathcal{M}_n\setminus\mathcal{PM}_n$. Furthermore, we
have
$$
W(\mathcal{EM }_n)-W(\mathcal{OM }_n)=W(\mathcal{PM}_n).
$$
\end{theorem}
\pf Given a matrix $A\in \mathcal{M}_n\setminus\mathcal{PM}_n $,
suppose that  $index(A)=i$.  We now have two cases. (1) We have
$csum_i(A)\geq 2$. (2) We have $1<i\leq dim(A)$, $csum_i(A)=1$,
$rsum_i(A)=0$, and $min_{i}(A)<min_{i-1}(A)$.

For Case (1), we obtain a new matrix $\Phi(A)$ from the matrix $A$
in the following way.  In A, replace the entry in row $min_i(A)$ of
column $i$ with zero. Then, insert a new zero row between row $i$
and row $i+1$ and insert a new column between column $i$ and $i+1$.
Let the new column be filled with all zeros except that the  entry
in row $min_i(A)$   is filled with $1$. In this case, we have
$\Phi(A)\in \mathcal{M}_n\setminus\mathcal{PM}_n $ with
$index(\Phi(A))=i+1$, $dim(\Phi(A))=dim(A)+1$ and
$rsum_1(\Phi(A))=rsum_1(A)$.

For Case (2), we may obtain a new matrix $\Phi(A)$ by reversing the
construction for Case  (1) as follows.  In $A$, replace the entry in
row $min_i(A)$ of column $i-1$ with  1. Then remove column $i$ and
row $i$. In this case, we have $\Phi(A)\in
\mathcal{M}_n\setminus\mathcal{PM}_n $ with $index(\Phi(A))=i-1$,
$dim(\Phi(A))=dim(A)-1$ and $rsum_1(\Phi(A))=rsum_1(A)$.

In both cases, the map $\Phi$ reverse the  parities   and preserve
the
  the weights of the matrices.  Hence, we obtain a desired
 parity reversing and weight preserving  involution on the set
$\mathcal{M}_n\setminus\mathcal{PM}_n$.   Note that if a matrix
$A\in
 \mathcal{M}_n$ is proper, then    there is exactly one $1$ in
 each column. Hence for each $A\in
 \mathcal{PM}_n$,  the parity of $A$ is even.
By applying  the involution, we can deduce that
$$ W(\mathcal{EM }_n)-W(\mathcal{OM }_n)=W(\mathcal{PM}_n).
$$

\qed

\begin{example}
Consider the following two matrices in $\mathcal{M}_6$:
$$
A=\begin{bmatrix}
1 &1   &0  & 0    \\
0 & 1   & 1 & 0  \\

0 & 0 & 1 &  0  \\
0 &  0 &0 & 1  \\
\end{bmatrix},
\,\,\, B= \begin{bmatrix}
1 &1   &0  &0& 0    \\
0 & 1   &0& 1 & 0  \\

0 & 0 & 1 &0&  0  \\
0 & 0 & 0 &0&  0  \\

0 &  0 &0 &0 &1  \\
\end{bmatrix}.
$$

For matrix $A$,  we have $index(A)=3$. Thus we have
$$
 \Phi(A)= \begin{bmatrix}
1 &1   &0  &\bf{0}& 0    \\
0 & 1   &0& \bf{1} & 0  \\

0 & 0 & 1 &\bf{0}&  0  \\
\bf{0} &\bf{ 0} & \bf{0} &\bf{ 0} &  \bf{0 } \\

0 &  0 &0 &\bf{0} &1  \\
\end{bmatrix},
$$
where the new inserted row and column are illustrated in bold.

For matrix $B$, we have $index(B)=4$.  Thus we have
$$
\Phi(B)=\begin{bmatrix}
1 &1   &0  & 0    \\
0 & 1   & 1 & 0  \\

0 & 0 & 1 &  0  \\
0 &  0 &0 & 1  \\
\end{bmatrix}.
$$
In fact, we have $\Phi(A)=B$ and $\Phi(B)=A$.
\end{example}

\section{Proof of  Conjecture \ref{conjecture}}
In this section, we will  show that the right-hand side of Formula
\ref{identity} is the generating function for the number of matrices
in $\mathcal{PM}_{n,k}$. Furthermore, we prove that there is a
bijection between the set $\mathcal{PM}_{n,k}$ and the set
$\mathcal{I}_{n,k}$, which implies Conjecture \ref{conjecture}.

Let
$$
A(t,z)= \sum_{n\geq 0}\prod_{i=1}^{n}(1-(1-t)^{i-1}(1-zt)).
$$
With the assumption that the empty product is as usual taken to be
$1$,    we have
$$
 A(t,z)=1+\sum_{n\geq 1}\prod_{i=1}^{n}
\sum_{j=1}^{i-1} ({i-1\choose j}+z{i-1\choose j-1})(-1)^{j-1}t^{j}.
$$
Define $A_n(z)$ to be the coefficient of $t^n$ in $A(t,z)$ for
$n\geq 1$, that is
\begin{equation}\label{A_n}
A(t,z)=1+\sum_{n\geq 1}A_n(z)t^n.\end{equation}
  Thus we have
$$
A_n(z)= \sum_{d=1}^{n}\sum_{ n_1+n_2+\ldots+n_d=n
}(-1)^{n-d}\prod_{j=1}^{d}({j-1\choose n_j}+z{j-1\choose n_j-1}),
$$
where the second summation is over all  compositions
$n_1+n_2+\ldots+n_d=n$ such that $n_j\geq 1$ for $j=1,2,\ldots, d$.
\begin{lemma}\label{A_n(z)}
For $n\geq 1$, we have
$$
A_n(z)=W(\mathcal{EM }_n)-W(\mathcal{OM }_n).
$$
\end{lemma}
\pf  Let $\mathcal{M}(n_1, n_2, \ldots, n_d)$ be the set of matrices
in $\mathcal{M}_n$ with $d$ columns   in which the column sum of
column $j$ is equal to $n_j$ for all $1\leq j\leq d$. In order to
get a matrix $A\in \mathcal{M}(n_1, n_2, \ldots, n_d)$, we should
choose $n_j$ places in column $j$ form $j$ places to arrange $1's$
for all $1\leq j\leq d$. we have two cases. (1) If $A_{1,j}=0$, then
we have ${j-1\choose n_j}$ ways to arrange $1$'s in column $j$. (2)
If $A_{1,j}=1$, then we have ${j-1\choose n_j-1}$ ways to arrange
the remaining $1$'s in column $j$. In the former case, column $j$
contributes $1$ to the weight of $A$. While in the latter case,
column $j$ contributes $z$ to the weight of $A$. Altogether, column
$j$ contributes ${j-1\choose n_j}+z{j-1\choose n_j-1}$ to the weight
of $\mathcal{M}(n_1, n_2, \ldots, n_d)$, which implies that $$
W(\mathcal{M}(n_1, n_2, \ldots, n_d))=\prod_{j=1}^{d}({j-1\choose
n_j}+z{j-1\choose n_j-1}).$$  It is clear that the parity of each
matrix in $\mathcal{M}(n_1, n_2, \ldots, n_d)$ is the parity of the
number $n-d$. When $d$ ranges from $1$ to $n$ and $n_1, n_2, \ldots,
n_d$ range over all compositions $n_1+n_2+\ldots+n_d=n$  such that
$n_j\geq 1$ for all $1\leq j\leq d$, we
   get the desired
result. \qed

Denote by $a_{n,k}$ the cardinality of the set $\mathcal{PM}_{n,k}$.
Assume that $a(0,0)=1$.

\begin{theorem}\label{A(t,z)}
 We have $$A(t,z)= \sum_{n\geq 0, k\geq
0}a_{n,k}z^kt^n=  \sum_{n\geq
0}\prod_{i=1}^{n}(1-(1-t)^{i-1}(1-zt)).$$
\end{theorem}
\pf  Combining Theorem  \ref{involution} and Lemma \ref{A_n(z)}, we
deduce that $A_n(z)=W(\mathcal{PM}_n)$ for $n\geq 1$.  Note that
$W(\mathcal{PM}_n)=\sum_{k=1}^{n}a_{n,k}z^k$  for $n\geq 1$.  Hence
we have  $$A(t,z)=1+\sum_{n\geq 1}A_n(z)t^n= \sum_{n\geq 0, k\geq
0}a_{n,k}z^kt^n, $$ which implies the desired result. \qed

From Theorem \ref{A(t,z)},   in order to prove Conjecture
\ref{conjecture}, it suffices to prove that $a_{n,k}=p_{n,k}$. In  a
matrix $A$, the operation of adding column $i$ to column  $j$ is
defined by increasing  $A_{k,j}$ by $A_{k,i}$  for each
$k=1,2,\ldots, dim(A)$. Note that a matrix $A\in \mathcal{M}_n$ is
proper if and only if it satisfies
\begin{itemize}
\item each column has exactly one $1$;
\item if $rsum_i(A)=0$, then we have $min_i(A)\geq min_{i-1}(A)$ for
$2\leq i\leq dim(A)$.
\end{itemize}
This observation will be  essential in the construction  of the
bijection between the set $\mathcal{PM}_{n,k}$ and the set
$\mathcal{I}_{n,k}$.

\begin{theorem}\label{bijection}
 There is  a bijection between the set $\mathcal{PM}_{n,k}$
and the set $\mathcal{I}_{n,k}$.
\end{theorem}

\pf Let $A$ be a matrix in the $\mathcal{PM}_{n,k}$, we can
construct  a matrix $A'$ in $\mathcal{I}_{n,k}$.  If there is no
zero rows in $A$, then we do nothing for $A$ and let $A'=A$. In this
case, the resulting matrix $A'$ is contained in $\mathcal{I}_{n,k}$.
Otherwise, we can construct a new upper triangular matrix  $A'$ by
the following {\em removal}  algorithm.
\begin{itemize}
\item Find the  least value $i$ such that row $i$ is a zero
 row. Then we obtain a new  upper triangular matrix by  adding column $i$
to  column  $i-1$ and remove column $i$ and row $i$.
\item Repeat the above procedure for the resulting matrix until
there is no zero row  in the resulting matrix.
\end{itemize}

 Clearly, the obtained matrix $A'$ is a matrix in
$\mathcal{I}_n$. Since the algorithm preserves the  sums of entries
in each non-zero rows of $A$, we have $rsum_1(A')=rsum_1(A)$. Hence,
the resulting matrix $A'$ is in $\mathcal{I}_{n,k}$.

Conversely, we can construct a matrix in $\mathcal{PM}_{n,k}$ from a
matrix in $\mathcal{I}_{n,k}$. Let $B$ be a matrix in the
$\mathcal{I}_{n,k}$. If the sum of entries in each  column is equal
to $1$, then we do nothing for $B$ and let $B'=B$. Otherwise, we can
construct a new upper triangular matrix  $B'$  by the following {\em
addition} algorithm.
\begin{itemize}
\item Find the  largest  value $i$ such that $csum_i(B)\geq 2$. Then we obtain a new  upper triangular matrix
by     decreasing    the entry in row $max_i(B)$ of column $i$ by
$1$, where $max_i(B)$ is defined to be the largest value $j$ such
that $B_{j,i}$ is non-zero.  Since $B$ is upper triangular, we have
$max_i(B)\leq i$.
\item Insert  one column between column $i$ and column $i+1$ and one zero
row between row $i$ and row $i+1$ such that the new inserted column
is filled with all zeros except that the entry on row $max_i(B)$ is
filled with $1$.
\item Repeat the above procedure for the resulting matrix until
there is no column  whose column sum is lager than $1$.
\end{itemize}
Clearly, the obtained matrix $B'$ is a matrix in $\mathcal{M}_n$.
From the construction of the above algorithm we know that the column
sum of each column  in $B'$ is equal to $1$. Furthermore, if row $j$
is a zero row, then we must have $min_j(B')\geq min_{j-1}(B')$.
Thus, the resulting matrix $B'$ is proper.
   Since the algorithm preserves the sums of entries in each
non-zero row of $B$,  we have   $rsum_1(B')=rsum_1(B)$. Hence, the
resulting matrix $B'$ is in $\mathcal{PM}_{n,k}$. This completes the
proof.

 \qed

\begin{example}
Consider a matrix $A\in \mathcal{PM}_{6,3}$. By applying the removal
algorithm, we get
$$
A=\begin{bmatrix}
1 &1   &\bf{0 } & 1 &0&0  \\
0 & 0   & \bf{1} & 0 &0&0 \\
{\bf0} & \bf{0}   & \bf{0} & \bf{0} &\bf{0}&\bf{0} \\
0 & 0   & \bf{0} & 0 &1&1 \\
0 & 0   & \bf{0} & 0 &0&0 \\
0 & 0   & \bf{0} & 0 &0&0 \\
\end{bmatrix}
\leftrightarrow  \begin{bmatrix}
1 &1   & 1 &\bf{0}&0  \\
0 & 1  & 0 &\bf{0}&0 \\

0 & 0  & 0 &\bf{1}&1 \\
\bf{0} & \bf{0}  & \bf{0 }&\bf{0}&\bf{0 }\\
0 & 0  & 0 &\bf{0}&0 \\
\end{bmatrix}
\leftrightarrow  \begin{bmatrix}
1 &1   & 1  &\bf{0}  \\
0 & 1  & 0  &\bf{0} \\

0 & 0  & 1  &\bf{1} \\

\bf{0} & \bf{0}  & \bf{0}  &\bf{0} \\
\end{bmatrix}
\leftrightarrow  A'=\begin{bmatrix}
1 &1   & 1     \\
0 & 1  & 0    \\
0 & 0  & 2    \\
\end{bmatrix},
$$
where the removed rows and columns are illustrated in bold at each
step of the removal algorithm. Conversely, given $A'\in
\mathcal{I}_{6,3}$, by applying addition algorithm, we can get $A\in
\mathcal{PM}_{6,3}$, where the inserted new rows and columns are
illustrated in bold at each step of the addition algorithm.
\end{example}

 Combining Theorems \ref{involution}, \ref{A(t,z)} and \ref{bijection},  we obtain
a combinatorial proof of Conjecture \ref{conjecture}. Note that
specializing $z=1$ implies a combinatorial proof of Formula
(\ref{pt}), which was proved by Bousquet-M$\acute{e}$lou et al.
\cite{melon} by using functional equations and the Kernel method.
 \vskip 5mm

\noindent{\bf Acknowledgments.}   The   author was supported by the
National Natural Science Foundation of China (no.10901141).


\end{document}